\documentclass{article}[12pt]
\usepackage{graphicx,pstricks,comment}
\usepackage[all]{xy}
\usepackage{subfigure}
\usepackage{rotating}
\usepackage{amsfonts}
\usepackage{amssymb}
\usepackage{amsmath}

\usepackage{indentfirst}
\usepackage{epsfig}
\usepackage{psfrag}
\usepackage{enumerate}
\usepackage{array}
\usepackage{theorem}
\usepackage{xcolor}

\newcommand{\C}{\mathbb{C}}
\newcommand{\R}{\mathbb{R}}
\newcommand{\N}{\mathbb{N}}
\newcommand{\Z}{\mathbb{Z}}

\newcommand{\Q}{\mathbb{Q}}
\newcommand{\quat}{\mathbb{H}}
\newcommand{\oct}{\mathbb{O}}
\newcommand{\K}{\mathbb{K}}

\newcommand{\HCn}{{\rm H}_\C^n}

\newcommand{\HR}{{\rm H}_{\R}}

\newcommand{\la}{\langle}
\newcommand{\ra}{\rangle}

\newtheorem{thm}{Theorem}

\newtheorem{lem}{Lemma}[section]
\newtheorem{prop}{Proposition}[section]
\newtheorem{dfn}{Definition}

\newtheorem{rmk}[thm]{Remark}

\newtheorem{question}{Question}

\newcounter{notes}%

\newcommand{\Pf}{{\em Proof}. }
\newcommand{\EPf}{\hfill$\Box$\vspace{.5cm}}
\numberwithin{equation}{section}

\title{Parabolic-preserving deformations of cusped hyperbolic lattices}

\author{Samuel A. Ballas, Julien Paupert, Pierre Will}

\begin{document}
\maketitle

\begin{abstract} We study deformations of non-cocompact lattices of ${\rm SO}(n,1)$ into ${\rm SU}(n,1)$ and ${\rm SO}(n+1,1)$. A necessary condition for these deformations to remain discrete and faithful (when $n \geqslant 3$) is for the parabolic subgroups to remain parabolic and discrete; we call such representations \emph{strongly parabolic-preserving}. We show that the figure-eight knot group admits a one-parameter family of Zariski-dense parabolic-preserving deformations into ${\rm SU}(3,1)$, with further deformations into ${\rm SU}(2,2)$. We also study the \emph{bending deformations} of the Bianchi groups (seen as subgroups of ${\rm SO}(3,1)$) along the modular surface into ${\rm SU}(3,1)$ and ${\rm SO}(4,1)$, and show that infinitely many of them are strongly parabolic-preserving in ${\rm SU}(3,1)$, while none are strongly parabolic-preserving in ${\rm SO}(4,1)$. Finally, for any $n \geqslant 3$, we show that there exist infinitely many non-commensurable cusped hyperbolic $n$-manifolds whose corresponding hyperbolic representation admits a 1-parameter family of parabolic-preserving deformations into ${\rm SU}(n,1)$. 
\end{abstract}

\section{Introduction}

This paper concerns an aspect of the deformation theory of discrete subgroups of Lie groups, namely that of non-cocompact (or \emph{cusped}) lattices in rank-one semisimple Lie groups. More specifically, we consider the following questions, given a discrete subgroup $\Gamma$ of a rank-one Lie group $H$:
\begin{itemize}
\item[(1)] Does $\Gamma$ admit any deformations in $H$?
\item[(2)] If so, do these deformations have any nice properties (e.g. remain discrete and faithful)?
\item[(3)] What if we replace $H$ with a larger Lie group $G$? 
\end{itemize}
Here we call  \emph{deformation} of $\Gamma$ in $H$ any continuous 1-parameter family of representations $\rho_t: \Gamma \longrightarrow H$ (for $t$ in some interval $(-\varepsilon,\varepsilon)$) satisfying $\rho_0=\iota$ (the inclusion of $\Gamma$ in $H$), and $\rho_t$ not conjugate to $\rho_{t'}$ for any $t \neq t' \in (-\varepsilon,\varepsilon)$. We say that $\Gamma$ is \emph{locally rigid} in $H$ if it does not admit any deformations into $H$.

When $H$ is a semisimple real Lie group without compact factors there are a variety of general local rigidity results which we now outline. Weil proved in \cite{W} that $\Gamma$ is locally rigid in $H$ if $H/\Gamma$ is compact and $H$ not locally isomorphic to ${\rm SL}(2,\R)$. Garland and Raghunathan extended this result to the case where $\Gamma$ is a non-cocompact lattice in a rank-one semisimple group $H$ not locally isomorphic to ${\rm SL}(2,\R)$ or ${\rm SL}(2,\C)$ (Theorem 7.2 of \cite{GR}). 

The exclusion of ${\rm SL}(2,\R)$ and ${\rm SL}(2,\C)$ is necessary. Generically, lattices in $H={\rm SL}(2,\R)$ admit many deformations in $H$. The identification of ${\rm PSL}(2,\R)$ with ${\rm Isom}^+(\HR^2)$ allow us to relate lattices in ${\rm SL}(2,\R)$ with hyperbolic structures which are in turn parameterized by the classical Teichm\"uller space when $\Gamma$ is a surface group. The case of \emph{quasi-Fuchsian} deformations of a discrete subgroup $\Gamma$ of  $H={\rm SL}(2,\R)$ into $G={\rm SL}(2,\C)$ is also classical, well-studied, and well understood by the Bers simultaneous uniformization theorem \cite{Bers}. In this setting, ${\rm PSL}(2,\C)$ can be identified with ${\rm Isom}^+(\HR^3)$ and the discrete group $\Gamma\subset {\rm SL}(2,\C)$ gives rise to a hyperbolic structure on the manifold $M\cong \Sigma\times \R$, where $\Sigma$ is a hyperbolic surface. Deforming $\Gamma$ in ${\rm SL}(2,\C)$ corresponds to deforming the hyperbolic structure on $M$. Such deformations are abundant and according the Bers simultaneous uniformization can be parameterized by a cartesian product of two copies of the classical Teichm\"uller space of the surface $\Sigma$. Notice that the existence of deformations into $G$ does not violate Weil's result as $G/\Gamma$ is not compact. 

This situation can be generalized to the case where $H={\rm SO}^0(n,1)\cong {\rm Isom}^+(\HR^n)$ and $G={\rm SO}^0(n+1,1)\cong {\rm Isom}^+(\HR^{n+1})$. In this setting, a lattice $\Gamma$ in $H$ gives rise to a hyperbolic structure on $M=\HR^n/\Gamma$. Again, regarding $\Gamma$ as a subgroup of ${\rm SO}^0(n+1,1)$ gives rise to a hyperbolic structure on $M\times \R$ and deformations of $\Gamma$ into $G$ correspond to deforming this hyperbolic structure. In this more general setting there is no general theorem that guarantees the existence of deformations of $\Gamma$ into $G$. However, when $n=3$ this deformation problem has been studied by Scannell \cite{Sc}, Bart--Scannell \cite{BSc}, and Kapovich \cite{Kap} who prove some rigidity results. 

Returning to the 3-dimensional case, many non-cocompact lattices $\Gamma$ in $H={\rm SL}(2,\C) \cong {\rm SO}^0(3,1)$ are known to admit deformations into $H$. In particular, when $\Gamma$ is torsion-free 
Thurston showed that for each cusp there exists a (real) 2-dimensional family of deformations of $\Gamma$ into $H$, called \emph{Dehn surgery deformations} (see Section 5.8 of \cite{Th}). Geometrically, in each of these families the commuting pair of parabolic isometries generating the corresponding cusp group is deformed to a pair of loxodromic isometries sharing a common axis. In particular these deformations are all non-discrete or non-faithful. If $\HR^3/\Gamma$ is an orbifold then the existence of deformations depends more subtly on the topology of the cusp cross-sections (see \cite{DM}).

Another case of interest in the context of deformations of geometric structures is that of \emph{projective deformations} of hyperbolic lattices, i.e. deformations of lattices $\Gamma$ of $H={\rm SO}(n,1)$ into $G={\rm SL}(n+1,\R)$. When $\Gamma$ is a torsion-free cocompact lattice in ${\rm SO}(n,1)$ such that the hyperbolic manifold $M={\rm H}^n_\R/\Gamma$ contains an embedded totally geodesic hypersurface $\Sigma$, Johnson and Millson showed in \cite{JM} that $\Gamma$ admits a 1-parameter family of deformations into ${\rm SL}(n+1,\R)$. They obtained these deformations, called \emph{bending deformations of $M$ along $\Sigma$}, by introducing an algebraic version of Thurston's bending deformations of a hyperbolic 3-manifold along a totally geodesic surface. 

This algebraic version is very versatile, and can be generalized in a variety of ways. For example, the hypothesis that $M$ is compact may be dropped, see \cite{BM}. Furthermore, the construction can be applied to the setting of other Lie groups and will provide us with a rich source of examples that are discussed in Section~\ref{bending}. In addition to deformations constructed via bending there are also instances of projective deformations that do not arise via the previously mentioned bending technique (see \cite{B1,B2,BDL}). On the other hand, despite the existence of these bending examples, empirical evidence compiled by Cooper--Long--Thistlethwaite \cite{CLT2} suggests that the existence of deformations into ${\rm SL}(4,\R)$ is quite rare for closed hyperbolic 3-manifolds.  

In another direction, \emph{complex hyperbolic quasi-Fuchsian} deformations of Fuchsian groups have also been extensively studied (see e.g. \cite{Sch}, the survey \cite{PP} and references therein).  With the above notation, this concerns deformations of discrete subgroups $\Gamma$ of $H$ into $G$, with $(H,G)=({\rm SO}(2,1),{\rm SU}(2,1))$ or $({\rm SU}(1,1),{\rm SU}(2,1))$. (Recall that the Lie groups ${\rm SL}(2,\R), {\rm SO}(2,1), {\rm SU}(1,1)$ are all isomorphic, up to index 2).  It turns out that for any $n \geqslant 2$, by work of Cooper--Long--Thistlethwaite \cite{CLT} there is an intricate relationship between projective deformations and complex hyperbolic deformations of finitely generated subgroups $\Gamma$ of ${\rm SO}(n,1)$, based on the fact that the Lie algebras of ${\rm SL}(n+1,\R)$ and ${\rm SU}(n,1)$ are isomorphic as modules over the ${\rm SO}(n,1)$ group ring. Specifically, they prove:

\begin{thm}[\cite{CLT}]\label{realprojcomphyp}
 Let $\Gamma$ be a finitely generated group, and let $\rho:\Gamma \longrightarrow {\rm SO}^0(n,1)$ be a smooth point of the representation variety ${\rm Hom}(\Gamma, {\rm SL}(n+1,\R))$. Then $\rho$ is also a smooth point of ${\rm Hom}(\Gamma, {\rm SU}(n,1))$, and near $\rho$ the real dimensions of ${\rm Hom}(\Gamma, {\rm SL}(n+1,\R))$ and ${\rm Hom}(\Gamma, {\rm SU}(n,1))$ are equal.
\end{thm}

The primary motivation for this article is to construct examples of complex hyperbolic deformations of real hyperbolic lattices $\Gamma$ that have nice algebraic and geometric properties. We are particularly interested in the case where the deformations remain discrete and faithful in a neighborhood of the original embedding. Indeed, from the point of view of geometric structures these correspond to deformations of the complex hyperbolic structure on the tangent bundle $TM$, where $M={\rm H}^n_\R/\Gamma$ is the hyperbolic manifold/orbifold corresponding to the lattice $\Gamma$ (see e.g. Theorem 3.1 of \cite{PT}).

We first recall that this property is automatically satisfied if $\Gamma$ is a cocompact lattice by the following result of Guichard, since a cocompact lattice is in particular convex-cocompact.

\begin{thm}[\cite{Gui}]\label{guichard}Let $G$ be a semisimple Lie group with finite center, $H$ a rank 1 subgroup of $G$, $\Gamma$ a finitely generated discrete subgroup of $H$ and denote $\iota : \Gamma \longrightarrow G$ the inclusion map. If $\Gamma$ is convex-cocompact then $\iota$ has a neighborhood in ${\rm Hom}(\Gamma,G)$ consisting entirely of discrete and faithful representations.
\end{thm}

On the other hand, when $\Gamma$ is a cusped lattice there is no reason to expect all deformations to locally remain discrete and faithful; typical examples of this are the Dehn surgery deformations of cusped 3-manifold groups into ${\rm SO}(3,1)$ mentioned earlier. As noted there, the failure of these deformations to remain both discrete and faithful is manifested geometrically by the fact that two commuting parabolic isometries in the initial lattice are deformed to a pair of commuting loxodromic isometries, which must therefore share the same axis, hence generate a non-discrete or non-faithful representation of $\Z^2$. Note that the same conclusion holds if the pair of commuting parabolic isometries is deformed to a pair of commuting elliptic isometries, showing the necessity of the following property for deformations of such a lattice to locally remain discrete and faithful, as long as the rank of the cusp group is at least two.

Recall that, from the Cartan--Killing classification of real semisimple Lie groups, any negatively curved symmetric space $X$ is a hyperbolic space ${\rm H}_\K^n$, with $\K=\R, \C, \quat$ or $\oct$ (and $n\geqslant 2$ if $\K=\R$, $n=2$ if $\K=\oct$). We refer the reader to \cite{CG} for general properties of these spaces and their isometry groups. In particular isometries of such spaces are roughly classified into the following 3 types: elliptic (having a fixed point in $X$), parabolic (having no fixed point in $X$ and exactly one on $\partial_\infty X$) or loxodromic (having no fixed point in $X$ and exactly two on $\partial_\infty X$). 
For our purposes we will need to distinguish between elliptic isometries with an isolated fixed point in $X$, which we call \emph{single-point elliptic} and elliptic isometries having fixed points on the boundary $\partial_\infty X$, which we call \emph{boundary elliptic}.

Given a discrete subgroup $\Gamma$ of ${\rm Isom}(X)$, we will call \emph{parabolic subgroup} of $\Gamma$ the stabilizer in $\Gamma$ of any boundary point fixed by a parabolic element of $\Gamma$. Note that such a subgroup may contain boundary elliptic elements; otherwise all its (non-identity) elements are parabolic. 

 \begin{dfn}\label{defpp} Let $X$ be a negatively curved symmetric space, $G={\rm Isom}(X)$, and $\Gamma$ a discrete subgroup of $G$. A representation $\rho:\Gamma \longrightarrow G$ is called \emph{parabolic-preserving} if $\rho(\gamma)$ is parabolic (resp. boundary elliptic) whenever $\gamma \in \Gamma$ is. We will call $\rho$ \emph{strongly parabolic-preserving} if it is parabolic-preserving and 
 $\rho(\Lambda)$ is discrete for any parabolic subgroup $\Lambda$ of $\Gamma$.
  \end{dfn}

{\bf Remark:} This is essentially what Bowditch called \emph{type-preserving} in \cite{Bo} (with more flexibility in the choice of parabolic subgroups). We change the name because ``type-preserving" has been widely used to mean that every parabolic (resp. elliptic, resp. loxodromic) isometry is mapped to a  parabolic (resp. elliptic, resp. loxodromic) isometry. The point of this weaker condition is that it is much easier to verify in practice, see the examples later in the paper.

\begin{question} Let $H \leqslant G$ be rank-one semisimple Lie groups and $\Gamma$ a cusped lattice in $H$. Is it true that any strongly parabolic-preserving deformation of $\Gamma$ into $G$ remains discrete and faithful in a neighborhood of the inclusion?
\end{question}

As far as we know this is an open question; we hope the following examples and non-examples will convince the reader that it is an interesting one.

\medskip

{\bf Remarks:}
\begin{itemize} 
\item The answer is ``yes" if we replace ``cusped" with ``cocompact" by Guichard's theorem above.  
\item As observed by Bowditch in \cite{Bo}, one needs to require ``strongly parabolic-preserving" rather than ``parabolic-preserving" in dimension at least three. In fact, it is noted there that Misha Kapovich observed that bending deformations of a cusped lattice of SO(3,1) into SO(4,1), along a \emph{cusped} totally geodesic embedded surface, are in some cases non-discrete even when restricted to the parabolic subgroup (see Proposition~\ref{CuspBendingNonSeparatingDim3} below for a precise statement). 

We give concrete examples of this for the Bianchi groups in Section~\ref{BianchiBending2}, as we find it interesting to contrast them with the case of bending the Bianchi groups from SO(3,1) into SU(3,1), where some of the bending deformations are in fact strongly parabolic-preserving (see Section~\ref{BianchiBending1}).
 
 \item On the other hand, it follows from Proposition 1.8 of \cite{Bo} that the answer to Question 1 is ``yes" when $H={\rm SO}(n,1)$ and $G={\rm SO}(n+1,1)$ for any $n \geqslant 2$. Indeed, in that case the peripheral subgroups (in the terminology of \cite{Bo}), that is the parabolic subgroups of the lattice $\Gamma < H$, are abelian groups of rank $n-1$. 
 
\item Bowditch also gives a concrete example of a sequence of representations $\rho_n$ of the free group $F_2 =\la a,b \ra$ into SO(4,1) converging to a representation $\rho_\infty$ such that: (1) $\rho_\infty(F_2)$ is discrete and contained in a conjugate of SO(2,1), (2) $\rho_n(F_2)$ is non-discrete for all $n$, and (3) $\rho_n(a)$ is parabolic for all $n$ as is $\rho_\infty(a)$; $\rho_n(b)$ is loxodromic for all $n$ as is $\rho_\infty(b)$. Thinking of the representations $\rho_n$ as belonging to a 1-parameter family of deformations of $\rho_\infty$, this gives an example of strongly parabolic-preserving deformations of $\rho_\infty(F_2)$ that are not both discrete and faithful in any neighborhood of $\rho_\infty(F_2)$. However this does not answer Question 1 in the negative, as $\rho_\infty(F_2)$ is not a lattice in SO(2,1), but an infinite-area discrete subgroup.

\end{itemize}

In this paper we produce deformations of cusped lattices of ${\rm SO}(n,1)$ into ${\rm SU}(n,1)$ and ${\rm SO}(n+1,1)$ for $n \geqslant 3$, and show that those into ${\rm SU}(n,1)$ are parabolic-preserving (and in some cases strongly parabolic-preserving). In the cases that we consider the deformations into ${\rm SO}(n+1,1)$ are not parabolic-preserving.

First, in Section~\ref{figure8} we consider deformations of the fundamental group of the simplest cusped hyperbolic 3-manifold, namely the figure-eight knot complement $S^3 \setminus K_8$. 
We denote $\Gamma_8=\pi_1(S^3 \setminus K_8)$, and $\rho_{\rm hyp}: \Gamma_8 \longrightarrow {\rm SO}(3,1)$ its \emph{hyperbolic representation}, i.e. the holonomy representation of the complete hyperbolic structure on $S^3 \setminus K_8$ (this is well-defined up to conjugation by Mostow rigidity). We obtain the following:

\begin{thm}\label{figure8main} Let $\Gamma_8$ be the figure-8 knot group and $\rho_{\rm hyp}: \Gamma_8 \longrightarrow {\rm SO}(3,1)$ its hyperbolic representation. Then there exists a 1-parameter family $\rho_u$ (with $u=e^{i\alpha} \in {\rm U}(1)$) of deformations of $\rho_{\rm hyp}=\rho_1$ into ${\rm SL}(4,\C)$, with image in (a conjugate of) ${\rm SU}(3,1)$ when $|\alpha| < 2\pi/3$, and ${\rm SU}(2,2)$ when $\alpha \in \pm (2\pi/3,\pi)$. The representations into ${\rm SU}(3,1)$ are parabolic-preserving, and $\rho_u$ has discrete image when $u=\pm 1, \pm i, e^{\pm i\pi/3}, e^{\pm 2i\pi/3}$.
\end{thm}
The proof is constructive in the sense that we give an explicit parametrization of the 1-parameter family of matrices in ${\rm SU}(3,1)$ generating the image of $\Gamma_8$. We hope that the explicit nature of these representations will allow the use of tools from \cite{GW} to prove discreteness and faithfulness of these representations in a neighborhood of $\rho_{\rm hyp}$.

In Section~\ref{bending} we review the algebraic description of \emph{bending deformations} of Johnson--Millson (\cite{JM}), and apply it to obtain bending deformations of any cusped hyperbolic $n$-manifold group $\Gamma < {\rm SO}(n,1)$ along a \emph{cusped} (properly embedded, totally geodesic) hypersurface into both ${\rm SU}(n,1)$ and ${\rm SO}(n+1,1)$. We determine under which conditions the deformations are parabolic and strongly parabolic. See Proposition~\ref{CuspBendingNonSeparatingDim3} below for a precise statement when $n=3$, and Proposition~\ref{CuspBendingNonSeparatingDimn} for general $n \geqslant 3$. In particular we obtain the following result:

\begin{thm}\label{df_all_dims}
For any $n\geqslant 3$ there exist infinitely many non-commensurable cusped hyperbolic $n$-manifolds whose corresponding hyperbolic representation admits a 1-parameter family of parabolic-preserving deformations into ${\rm SU}(n,1)$. 
\end{thm} 

Here, two groups $\Gamma,\Gamma'\subset H$ are \emph{commensurable} (in the wide sense) if $\Gamma\cap g\Gamma'g^{-1}$ has finite index in both $\Gamma$ and $g\Gamma'g^{-1}$ for some $g  \in H$. (The incommensurability conclusion ensures that in each dimension $n$ the manifolds in Theorem \ref{df_all_dims} are quite distinct in the sense that they are not obtained by taking covering spaces of a single example). The manifolds in question are obtained by taking appropriate  covers of arithmetic orbifolds of \emph{simplest type}, which by deep results of Bergeron--Haglund--Wise (\cite{BHW}) and McReynolds--Reid--Stover (\cite{MRS}) have only torus cusps and contain a cusped embedded totally geodesic hypersurface which intersects each cusp cross-section in at most one component (see Section~\ref{BendingAllDims}). We then carefully analyze the bending deformations along such a hypersurface to show that they are parabolic-preserving. 

Finally, we explicitly compute these bending deformations for the Bianchi groups into ${\rm SU}(3,1)$ and  ${\rm SO}(4,1)$, where the cusped surface along which we bend is the modular surface. As above, we  hope that the explicit nature of these representations will allow the use of tools from \cite{GW} to prove discreteness and faithfulness of these representations in a neighborhood of $\rho_{\rm hyp}$. 

The following  summarizes our results about these deformations; see Propositions~\ref{BianchiDefsSU31} and \ref{BianchiDefsSO41} for more precise statements. The Bianchi groups are defined, for $d\geqslant 1$ squarefree, by ${\rm Bi}(d)={\rm PSL}(2,\mathcal{O}_d) < {\rm PSL}(2,\C)$, where $\mathcal{O}_d$ denotes as usual the ring of integers of $\Q (i\sqrt{d})$. We regard these as subgroups of ${\rm SO}(3.1)$ (using an explicit isomorphism ${\rm PSL}(2,\C) \simeq {\rm SO}^0(3.1)$). It is well known that the modular surface ${\rm H}^2_\R/{\rm PSL}(2,\Z)$ is embedded in the Bianchi orbifold ${\rm H}^3_\R/{\rm Bi}(d)$ for all $d \neq 1,3$. (The deformations of ${\rm Bi}(d)$ into ${\rm SU}(3,1)$ when $d=1,3$ were completely classified in \cite{PT}).

\begin{thm}\label{Bianchi} For all $d \neq 1,3$, the bending deformations of ${\rm Bi}(d) < {\rm SO}(3,1)$ into ${\rm SU}(3,1)$ along the modular surface are parabolic-preserving; they are moreover strongly parabolic-preserving whenever $d$ is congruent to 1 or 2  modulo  4.
\end{thm}
 
 {\bf Acknowledgements:} The second author would like to thank Teddy Weisman, Jeff Danciger and Misha Kapovich for helpful conversations concerning the material in this paper.
 
\section{Background on real and complex hyperbolic space}

We start with a brief review of definitions and basic facts about the geometry of real and complex hyperbolic spaces, their isometries and the structure of their boundary at infinity. In fact we will center the discussion on complex hyperbolic space, and observe that this contains real hyperbolic space by restricting the vector space $\C^{n+1}$ to $\R^{n+1}$ (or more precisely, any $\R$-linear $n$-dimensional subspace on which the Hermitian form takes real values).
See \cite{CG} or \cite{Go} for a more detailed discussion.

 \medskip
 
{\bf Projective models:} Consider the vector space $\C^{n+1}$ endowed with a Hermitian form $\langle \cdot \, , \cdot \rangle$ of signature $(n,1)$.
	Let $V^-=\left\lbrace Z \in \C^{n,1} | \langle Z , Z \rangle <0 \right\rbrace$ and denote $\pi: \C^{n+1}-\{0\} \longrightarrow \C{\rm P}^n$ the projectivization map.
	Define ${\rm H}_\C^n$ to be $\pi(V^-) \subset \C{\rm P}^n$, endowed with the distance $d$ given by:
	﻿
	\begin{equation}\label{dist}
		\cosh ^2 \frac{1}{2}d(\pi(X),\pi(Y)) = \frac{|\langle X, Y \rangle|^2}{\langle X, X \rangle  \langle Y, Y \rangle}
	\end{equation}
	﻿
	This distance is in fact induced by a complete Riemannian metric (called the \emph{Bergman metric}) which has negative sectional curvatures $1/4$-pinched (between $-1$ and $-1/4$ in the above normalization).
	
	Different choices of Hermitian forms of signature $(n,1)$ give rise to different models of $\HCn$. For the computations in Section~\ref{bending} we will use the  \emph{Siegel model}, where the Hermitian form is given by $\la Z,Z\ra=2{\rm Re} (z_1\overline{z_{n+1}}) +|z_2|^2+\cdots +|z_n|^2$. This form provides especially nice coordinates for the stabilizer of a particular point at infinity, in the form of upper-triangular matrices. (The Hermitian form that we will use in Section~\ref{figure8} to study deformations of the figure-eight knot group will not be as nice and will in fact vary with the parameter).
	
﻿\medskip
	
{\bf Isometries:} From \eqref{dist} it is clear that ${\rm PU}(n,1)$ acts by isometries on ${\rm H}_\C^n$, where ${\rm U}(n,1)$ denotes the subgroup of ${\rm GL}(n+1,\C)$ preserving $\langle \cdot , \cdot \rangle$, and ${\rm PU}(n,1)$ its image in ${\rm PGL}(n+1,\C)$. In fact, PU($n$,1) is the group of holomorphic isometries of ${\rm H}_\C^n$ (the full group of isometries ${\rm Isom}({\rm H}_\C^n)$ is ${\rm PU}(n,1) \rtimes \Z/2$, where the $\Z/2$ factor corresponds to a \emph{real reflection}, that is, an antiholomorphic involution). 

Holomorphic isometries of $\HCn$ are of three types, depending on the number and location of their fixed points. Namely, $g \in {\rm PU}(n,1) - \{ {\rm Id} \}$ is :
\begin{itemize}
\item \emph{elliptic} if it has a fixed point in ${\rm H}_\C^n$
		
\item \emph{parabolic} if it has (no fixed point in ${\rm H}_\C^n$ and) exactly one fixed point in $\partial{\rm H}_\C^n$
		
\item \emph{loxodromic}: if it has (no fixed point in ${\rm H}_\C^n$ and) exactly two fixed points in $\partial{\rm H}_\C^n$
\end{itemize}
	﻿
A parabolic isometry is called {\it unipotent} (or a \emph{Heisenberg translation}) if it has a unipotent lift in ${\rm U}(n,1)$; if not it is called {\it ellipto-parabolic}.  In dimensions $n>1$, unipotent isometries are either  {\it 2-step} or {\it 3-step}, according to whether the minimal polynomial of their unipotent lift is $(X-1)^2$ or $(X-1)^3$ (see section 3.4 of \cite{CG}). We will usually call 2-step nilpotent isometries \emph{vertical Heisenberg translations} and 3-step nilpotent isometries \emph{horizontal Heisenberg translations}. This terminology reflects the position of the subspaces preserved by each type of isometry in the Heisenberg group, relative to its standard contact structure. Any two horizontal (resp. vertical) Heisenberg translations are conjugate under an isometry.

The following criterion (contained in Theorem 3.4.1 of \cite{CG}) is very convenient to understand isometries given in arbitrary matrix form:

\begin{prop}\label{IsomClasses} Let $A \in {\rm U}(n,1) \setminus \{{\rm Id}\}$ and $g \in {\rm PU}(n,1)$ the corresponding isometry of $\HCn$.
\begin{itemize}
\item $g$ is elliptic $\iff$ $A$ is diagonalizable with all eigenvalues of unit norm.
\item $g$ is loxodromic $\iff$ $A$ is diagonalizable with exactly $n-1$ eigenvalues of unit norm.
\item $g$ is parabolic $\iff$ $A$ is not diagonalizable. All eigenvalues of $A$ then have unit norm.
\end{itemize}
\end{prop}

\medskip
	
{\bf The Siegel model:}
As mentioned above, this model corresponds to the Hermitian form given by the matrix:
	\begin{equation}\label{Hermitian}
		H=
		\begin{bmatrix}
			0 & 0 & 1 \\
			0 & {\rm Id} & 0 \\
			1 & 0 & 0 
		\end{bmatrix}
	\end{equation}
	In this model, any point  $p \in \HCn$ admits a unique lift to $\C^{n,1}$ of the following form, called its \emph{standard lift}:
	\begin{equation}
		{\bf p}=\begin{bmatrix}
			(-|Z|^2-u+it)/2\\Z\\1
		\end{bmatrix}\mbox{ with } (Z,t,u)\in\C^{n-1} \times \R \times (0,\infty) .
	\end{equation}
The coordinates $(Z,t,u)$ are called \textit{horospherical coordinates} of $p$. The boundary at infinity $\partial_\infty\HCn$ is the level set $\lbrace u=0\rbrace$, together with the distinguished point at infinity, given by
	﻿$$
	p_\infty\sim\left[
	\begin{array}{c}
		1 \\ 0 \\ 0
	\end{array}\right].
	$$
Level sets $\lbrace u=u_0\rbrace$ with fixed $u_0>0$ are called \emph{horospheres based at $p_\infty$}. The punctured boundary at infinity $\partial_\infty\HCn \setminus \{ p_\infty \}$ is a copy of the \emph{generalized Heisenberg group}, identified with $\C^{n-1} \times \R$ with group law given by:
	﻿
	\begin{equation}
		\label{Heisprod}(Z_1,t_1)\cdot (Z_2,t_2)=(Z_1+Z_2,t_1+t_2+2{\rm Im}(Z_1 Z_2^*)).
	\end{equation}
(Note that for the obvious real hyperbolic subspace obtained by taking $Z \in \R^{n-1}$ and $t=0$ this reduces to $\R^{n-1}$ with its usual Euclidean structure). 

The stabilizer of $p_\infty$ in ${\rm PU}(n,1)$ consists of upper-triangular matrices, and is generated by the following 3 types of isometries: Heisenberg translations $T_{(Z,t)}$ ($(Z,t)\in \C^{n-1} \times \R$), Heisenberg rotations $R_U$ ($U \in {\rm U}(n-1)$) and Heisenberg dilations $D_r$ ($r >0$), where:
	﻿
	\begin{equation}\label{stabinf}
		\begin{array}{ccc}
			T_{(z,t)}=
			\begin{bmatrix}
				1 & -Z^* & -(|Z|^2-it)/2 \\
				0 & {\rm Id} & Z \\
				0 & 0 & 1
			\end{bmatrix}
			&
			R_U=
			\begin{bmatrix}
				1 & 0 & 0 \\
				0 & U & 0 \\
				0 & 0 & 1
			\end{bmatrix}
			&
			D_r=
			\begin{bmatrix}
				r & 0 & 0 \\
				0 & {\rm Id} & 0 \\
				0 & 0 & 1/r
			\end{bmatrix}.
		\end{array}
	\end{equation}
In Heisenberg coordinates, these correspond to the following:
\begin{itemize}
\item $T_{(Z,t)}$ is left multiplication by $(Z,t)$: $(W,s)\longmapsto (Z,t)\cdot (W,s)$,
\item $R_U$ is given by $(W,s)\longmapsto (UW,s)$,
\item $D_r$ is the Heisenberg dilation $(W,s)\longmapsto (rW,r^2s)$.
\end{itemize}
Heisenberg translations and rotations preserve each horosphere based at $p_\infty$ whereas Heisenberg dilations permute horospheres based at $p_\infty$. 

\section{Deformations of the figure-8 knot group into ${\rm SU}(3,1)$}\label{figure8}

In this section we construct a family of parabolic-preserving deformations of the hyperbolic representation of the figure-8 knot group into ${\rm SU}(3,1)$, and prove Theorem~\ref{figure8main} from the introduction. 

Consider $\Gamma_8=\pi_1(S^3 \setminus K_8)$ where $K_8$ is the figure-8 knot, and denote $\rho_{\rm hyp}: \Gamma_8 \longrightarrow {\rm SO}(3,1)$ the holonomy of the complete hyperbolic structure on $S^3 \setminus K_8$.
Recall that in the presence of a smoothness hypothesis on the relevant representation varieties, Theorem \ref{realprojcomphyp} implies that the existence of deformations of $\rho_{\rm hyp}$ into ${\rm SL}(4,\R)$ guarantees the existence of deformations of $\rho_{\rm hyp}$ into ${\rm SU}(3,1)$. Work of Ballas--Danciger--Lee \cite{BDL} shows that the smoothness hypothesis is guaranteed in the presence of a cohomological condition. Specifically, they prove the following.

\begin{thm}[\cite{BDL}]\label{infrigimpliessmooth} Let $M$ be an orientable complete finite volume hyperbolic manifold with fundamental group $\Gamma$, and let $\rho_{\rm hyp}:\Gamma \longrightarrow {\rm SO}(3,1)$ be the holonomy representation of the complete hyperbolic structure. If $M$ is infinitesimally projectively rigid rel boundary, then $\rho_{\rm hyp}$ is a smooth point of ${\rm Hom}(\Gamma, {\rm SL}(4,\R))$ and its conjugacy class is a smooth point of $\chi (\Gamma, {\rm SL}(4,\R))$.
\end{thm}

Roughly speaking,  \emph{infinitesimally projectively rigid rel boundary} is a cohomological condition that says that a certain induced map from the twisted cohomology of $M$ into the twisted cohomology of $\partial M$ is an injection. For a more precise definition, see \cite{HP}. By work of Heusener--Porti \cite{HP}, it is known that the figure-8 knot complement is infinitesimally rigid rel boundary, and so we can apply Theorems \ref{infrigimpliessmooth} and \ref{realprojcomphyp} to produce deformations of $\rho_{\rm hyp}$ into ${\rm SU}(3,1)$. However, there is no reason why these representations should be parabolic-preserving, and in many cases the deformations will not have this property.

Fortunately, work of the first author (see \cite{B1,B2}) provides a family of deformations of $\rho_{\rm hyp}$ into ${\rm SL}(4,\R)$ whose corresponding deformations into ${\rm SU}(3,1)$ are parabolic preserving.   

\begin{thm}[\cite{B1},\cite{B2}] Let $\Gamma_8$ be the figure-8 knot group. Then there exists a 1-parameter family of discrete, faithful deformations of $\rho_{\rm hyp}$ into ${\rm SL}(4,\R)$. 
\end{thm}

The construction of this 1-parameter family can be found in \cite{B1} and ultimately constructs a curve $\rho_t$ of representations of $\Gamma_8$ into ${\rm SL}(4,\R)$ containing the hyperbolic representation $\rho_{\rm hyp}$ at $t=1/2$. In fact, allowing the parameter $t$ to take complex values gives a 1-complex parameter family of representations into ${\rm SL}(4,\C)$. Moreover, it turns out that taking $2t$ to be a unit complex number $u$ gives a 1-parameter family of representations into ${\rm SU}(3,1)$. (The reason for this choice of value of the parameter is that the eigenvalues of one of the peripheral elements in $\rho_t(\Gamma_8)$ are 1 and a power of $2t$, see Section 6 of \cite{B2}). 

We now give explicit matrices for the generators and Hermitian form for this family, using the presentation and notation of Section 6 of \cite{B2}. There, the following presentation of $\Gamma_8$ was used:

\begin{equation}\label{mnpres}
\Gamma_8 = \la m,n \, | \, mw=wn\ra, \ {\rm where} \ w=[n,m^{-1}].
\end{equation}

The family of representations $\rho_u : \Gamma_8 \longrightarrow {\rm SL}(4,\C)$ is defined by $\rho_u(m)=M_u$ and $\rho_u(n)=N_u$, where:

\begin{equation}\label{generators}
\begin{array}{lcr}
M_u=\left(\begin{array}{cccc}
1 & 0 & 1 & u/2-1 \\
0 & 1 & 1 & u/2 \\
0 & 0 & 1 & (u+1)/2 \\
0 & 0 & 0 & 1
\end{array}\right)
&
{\rm and}
&
N_u=\left(\begin{array}{cccc}
1 & 0 & 0 & 0 \\
2(1+\bar{u}) & 1 & 0 & 0 \\
2 & 1 & 1 & 0 \\
1 & 1 & 0 & 1 
\end{array}\right)
\end{array}
\end{equation}   
 
When $|u|=1$, the group $\rho_u(\Gamma_8)$ preserves the Hermitian form $H_u$ on $\C^4$ given by $H(X,Y)=X^TJ_u\bar{Y}$, where:

\begin{equation}\label{form}
J_u=\left(\begin{array}{cccc}
1+(u+\bar{u})/2 & -1-(u+\bar{u})/2 & 1+u & -3 -2(u+\bar{u})-\bar{u}^2 \\
-1-(u+\bar{u})/2 & 1+(u+\bar{u})/2 & -1 - u & 1+u \\
1+\bar{u} & -1-\bar{u} & 4+2(u+\bar{u}) & -4-2(u+\bar{u}) \\
-3-2(u+\bar{u})-u^2 & 1+\bar{u} & -4-2(u+\bar{u}) & 4+2(u+\bar{u})
\end{array}\right)
\end{equation}

\begin{lem} The form $H_u$ has signature (3,1) for all $u=e^{i\alpha}$ with $|\alpha| < 2\pi/3$, and signature (2,2) when $\alpha \in \pm (2\pi/3,\pi)$.
\end{lem}

\Pf Computing the determinant of $J_u$ gives: 
$$
\begin{array}{rcl}
{\rm det} \, J_u & = &-96-83(u+\bar{u})-53(u^2+\bar{u}^2)-24(u^3+\bar{u}^3)-7(u^4+\bar{u}^4)-(u^5+\bar{u}^5) \\
 & = & -96-166{\rm cos}(\alpha)-106{\rm cos}(2\alpha)-48{\rm cos}(3\alpha)-14{\rm cos}(4\alpha) -2{\rm cos}(5\alpha)\\
 & = & -4 ({\rm cos}(\alpha)+1)^2(2{\rm cos}(\alpha)+1)^3.
\end{array}
$$
The latter function of $\alpha$ is negative for $|\alpha| < 2\pi/3$, and positive for $\alpha \in \pm (2\pi/3,\pi)$; the result then follows by noting that $H_u$ has signature (3,1) when $u=1$ (corresponding to the hyperbolic representation), and (2,2) for e.g. $u=\pm 3\pi/4$. \EPf

\begin{lem} The representations $\rho_u$ are pairwise non-conjugate in ${\rm SL}(4,\C)$.
\end{lem}
\Pf A straightforward computation gives: ${\rm Tr} \, M_uN_u=6+u$. \EPf

\begin{lem} For any $u=e^{i\alpha}$ with $|\alpha| < 2\pi/3$, the representation $\rho_u : \Gamma_8 \longrightarrow {\rm SU}(3,1)$ is parabolic-preserving.
\end{lem}

\Pf The peripheral subgroup $\Gamma_\infty$ of $\Gamma_8$ is generated by $m$ and $l=ww^{op}=nm^{-1}n^{-1}m^2n^{-1}m^{-1}n$, with the notation of the presentation~(\ref{mnpres}) (see \cite{B1}). Now $M_u=\rho_u(m)$ is unipotent for all $u$, and a straightforward computation (using eg Maple) shows that $L_u=\rho_u(l)$ is non-diagonalizable (with eigenvalues $(u,u,u,\bar{u}^3)$) for all $u$, hence parabolic by Proposition~\ref{IsomClasses}. Since $\Gamma_\infty \simeq \Z ^2$, all elements of $\rho_u(\Gamma_\infty)=\langle \rho_u(m),\rho_u(l) \rangle$ will also remain parabolic for $u$ in a neighborhood of 1. \EPf

For future reference the explicit matrix form of $L_u=\rho_u(l)$ is the following (using $\bar{u}=u^{-1}$ to save space):

$$L_u=\left( \begin{array}{cccc}
(u-1-\bar{u}-\bar{u}^2)/2 & (u+1+\bar{u}+\bar{u}^2)/2 & -1-\bar{u}^2 & (5u+6+2\bar{u}+3\bar{u}^2)/2 \\
(u-1-\bar{u}-\bar{u}^2-2\bar{u}^3)/2 & (u+1+\bar{u}+\bar{u}^2+2\bar{u}^3)/2 & 1-2\bar{u}+\bar{u}^2-2\bar{u}^3 & (7u+4+10\bar{u}+\bar{u}^2+6\bar{u}^3)/2 \\
0 & 0 & u & 0 \\
0 & 0 & 0 & u
\end{array} \right)
$$

The representations $\rho_u$ don't quite have all their matrix entries in $\Z[u,u^{-1}]$, however they still satisfy the following integrality property which ensures discreteness of some of their images:

\begin{prop}\label{figure8prop} 
\begin{enumerate}
\item For any $u \in \C$, ${\rm Tr} \, \rho_u(\Gamma_8) \subseteq \Z[u,u^{-1}]$.
\item If $u$ is a unit in a number field $k$, then $\rho_u(\Gamma_8)$ has a finite-index subgroup contained in ${\rm SL}(4,\mathcal{O}_k)$.
\item $\rho_u(\Gamma_8)$ is a discrete subgroup of ${\rm SL}(4,\C)$ when $u=\pm 1, \pm i, e^{\pm i\pi/3}, e^{\pm 2i\pi/3}$.
\end{enumerate}
\end{prop}
\Pf \begin{enumerate}
\item Since $\rho_u(\Gamma_8) < {\rm SL}(4,\C)$ is generated by the two matrices $M_u$ and $N_u$, it suffices to check that the traces of a finite number of words in these generators lie in $\Z[u,u^{-1}]$, namely the 17 words listed in Theorem 4.2 of \cite{GL} or Proposition 5.2 of \cite{Te}. For concreteness we list the first few of these words in $m,n=M_u,N_u$: $m, n, mn, mn^{-1}, mnm, mn^{-1}m, [m,n], mnmn^{-1}, mn^{-1}mn$, with corresponding traces $4,4,6+u, 3, 9+3u, 3+u^{-1}, 3, 6+u^{-1}, 6+u^{-1}$. The remaining words are analogous.
\item This follows immediately from the previous item and Theorem~2.4 of \cite{BL}, since if $u$ is a unit in $k$ then $\Z[u,u^{-1}] \subseteq \mathcal{O}_k$.
\item The result follows from the previous item and discreteness of the ring $\mathcal{O}_k$ in $\C$ when $k=\Q, \Q(i), \Q(e^{i\pi/3})$.
\end{enumerate} \EPf

\section{Bending deformations}\label{bending}

\subsection{The general construction}

Let $\Gamma$ be a lattice in ${\rm SO}(n,1)$; then $M=\HR^n/\Gamma$ is a finite-volume hyperbolic $n$-orbifold. Suppose that $M$ contains an embedded connected totally geodesic hypersurface $\Sigma=\HR^{n-1}/\Delta$ (with $\Delta = {\rm Stab}_\Gamma(\Sigma)$, a lattice in a conjugate of ${\rm SO}(n-1,1)$). Johnson--Millson gave in \cite{JM} an algebraic description of a one-parameter family of representations $\rho_t: \Gamma \longrightarrow G$, ($t \in \R$, with $\rho_0=\rho_{hyp}$), called \emph{bending deformations of $M$ along $\Sigma$} (or of $\Gamma$ along $\Delta$), under the assumption that $G$ is a Lie group containing ${\rm SO}(n,1)$ such that the centralizer of ${\rm SO}(n-1,1)$ in $G$ has dimension 1. When $G={\rm SL}(n,\R)$, ${\rm SU}(n,1)$ or ${\rm SO}(n+1,1)$ and $M$ a manifold (that is, $\Gamma$ torsion-free), they proved that these representations are indeed deformations of $\Gamma$, that is that they are pairwise non-conjugate in $G$. 

We now review this construction. The hypersurface $\Sigma$ provides a decomposition of $\Gamma$ into either an amalgamated free product or an HNN extension, depending on whether or not $\Sigma$ is separating. Using this decomposition we can construct a family $\rho_t:\Gamma \to G$ such that $\rho_0=\iota$, where $\iota$ is the inclusion of $\Gamma$ into $G$, as follows. We start by recalling the definitions of amalgamated free product and HNN extension:

$$
(1) \ \Gamma=\Gamma_1 *_\Delta \Gamma_2 \ \ \ \ {\rm or } \ \ \ \ (2) \ \Gamma=\Gamma_1 *_z 
$$
In the first case we start with presentations $\Gamma_1 = \langle S_1 \, | \, R_1 \rangle$ and $\Gamma_2 = \langle S_2 \, | \, R_2 \rangle$ and two homomorphisms $\varphi_1: \Delta \longrightarrow \Gamma_1$ and $\varphi_2: \Delta \longrightarrow \Gamma_2$. (Typically, these can be taken as the inclusion maps of a common subgroup $\Delta$). The amalgamated product $\Gamma$ admits the following presentation:

$$ (1) \ \Gamma = \Gamma_1 *_\Delta \Gamma_2 = \langle S_1,S_2 \, | \, R_1, R_2, \varphi_1(\gamma)=\varphi_2(\gamma) \, (\forall \gamma \in \Delta) \rangle.
$$

In the second case we start with a group $\Gamma_1$ with presentation $\Gamma_1 = \langle S_1 \, | \, R_1 \rangle$ and assume that $\Gamma_1$ contains two isomorphic subgroups $\Delta,\Delta'$ with an isomorphism $\varphi: \Delta \longrightarrow \Delta'$. The corresponding HNN extension is then defined, for some generator $z \notin \Gamma_1$, by:

$$(2) \ \Gamma = \Gamma_1 *_z = \langle S_1,z \, | \, R_1, z^{-1}\gamma z=\varphi(\gamma) \, (\forall \gamma \in \Delta) \rangle.
$$

Now denote $Z_G(\Delta) = \{g_t \}_{t \in \R}$ (a smooth parametrization of the connected component of the identity of) the centralizer in $G$ of $\Delta$. (Since $\Delta$ is a lattice in ${\rm SO}(n-1,1)$, the latter is the same as the centralizer of ${\rm SO}(n-1,1)$ in $G$).

If $\Sigma$ is separating, then $M\backslash \Sigma$ consists of two connected components $M_1$ and $M_2$, with fundamental groups $\Gamma_1$ and $\Gamma_2$ respectively. In this case $\Gamma=\Gamma_1\ast_{\Delta}\Gamma_2$.  
The group $\Gamma$ is generated by $\Gamma_1\cup\Gamma_2$ and we define 

\begin{equation}\label{BendingAmalgam}
\rho_t(\gamma)=\left\{\begin{matrix}
                       \iota(\gamma) & \gamma\in \Gamma_1\\
                       g_t\iota(\gamma)g_t^{-1} & \gamma\in \Gamma_2
                      \end{matrix}\right\}
\end{equation}

on this generating set. Since $g_t$ centralizes $\Delta$ we see that the relations coming from the amalgamated product decomposition are satisfied, and so $\rho_t:\Gamma\to G$ is well defined. 

If $\Sigma$ is non-separating, then $M'=M\backslash \Sigma$ is connected. If we let $\Gamma'$ be the fundamental group of $M'$ then we can arrive at the decomposition $\Gamma=\Gamma'\ast_z$. In this case $\Gamma$ is generated by $\Gamma'\cup\{z\}$, where $z$ is a free letter and we define $\rho_t$ on generators as

\begin{equation}\label{BendingHNN}
\rho_t(\gamma)=\left\{\begin{matrix}
                              \iota(\gamma) & \gamma\in \Gamma'\\
                              g_t \iota(\gamma) & \gamma=z
                             \end{matrix}\right\}
\end{equation}

Again, since $g_t$ centralizes $\Delta$ we see that the relations for the HNN extension are satisfied and so $\rho_t:\Gamma\to G$ is well defined. These representations are called \emph{bending deformations of $\Gamma$ along $\Delta$} (or of $M$ along $\Sigma$),

\subsection{Bending hyperbolic $n$-orbifolds along a totally geodesic hypersurface with torus cusps}

We now examine in more detail how a torus cusp group deforms under bending deformations along a hypersurface containing the cusp. 
For simplicity, we first state and prove the result in dimension 3. We warm up with the following preliminary result.

\begin{lem}\label{RtimesS1} Let $T=(a,0)$ and $U=(b,\theta)$ in $\R \times S^1$, with $a,b\neq 0$ and $S^1$ identified with $\R/2\pi\Z$.
\begin{enumerate}
\item If $a/b \notin \Q$ then $\la T,U \ra$ is non-discrete in $\R \times S^1$ and isomorphic to $\Z^2$.
\item If $a/b \in \Q$ and $\theta \in \pi \Q$ then $\la T,U \ra$ is discrete in $\R \times S^1$ and not isomorphic to $\Z^2$.
\item If $a/b \in \Q$ and $\theta \notin \pi \Q$ then $\la T,U \ra$ is non-discrete in $\R \times S^1$ and isomorphic to $\Z^2$.
\end{enumerate}
\end{lem}

\Pf \begin{enumerate}
\item Since $S^1$ is compact, if $\Gamma$ is discrete in $\R \times S^1$ then $p_1(\Gamma)$ is discrete in $\R$, denoting $p_1$ projection to the first factor. The statement about non-discreteness follows, since $\la a,b \ra$ is non-discrete (in fact, dense) in $\R$ when $a/b \notin \Q$. Moreover in that case $\la a,b \ra$ is isomorphic to $\Z^2$, hence so is $\la T,U \ra$. 
\item Since $\theta \in \pi \Q$, some power $U^k$ is of the form $(kb,0)$. Then $\la T,U^k \ra$ is discrete in $\R \times \{ 0\}$ (as $a/b \in \Q$); since this has finite index in $\la T, U \ra$, the latter is also discrete in $\R \times S^1$. Note that $T $ and $U^k$ (hence $T$ and $U$) satisfy some non-trivial relation in addition to commutation.
\item Since $a/b \in \Q$, after possibly replacing $T,U$ by powers and rescaling we may assume that $a=b=1$. Then, for any $m,n \in \Z$, $T^m U^n=(m+n,n\theta)$. The latter is a dense subset of $\Z \times S^1$ (since at each constant height $m+n=C$ we have a dense subset $T^{C-n} U^n = \{(C, n \theta ) \, | \, n \in \Z\}$ of $\{ C \} \times S^1$), hence is not discrete in $\R \times S^1$.
\end{enumerate}\EPf

\begin{prop}\label{CuspBendingNonSeparatingDim3} Let $M=\HR^3/\Gamma$ be a finite-volume hyperbolic 3-orbifold with a torus cusp represented by $p_\infty \in \partial_\infty \HR^3$, and $\Sigma=\HR^2/\Delta$ a properly embedded, non-separating totally geodesic surface in $M$ containing (the image of) $p_\infty$. 
Let $\rho_t: \Gamma \longrightarrow G$ ($t \in \R$) denote the bending deformations of $\Gamma$ along $\Delta$, with $G={\rm SU}(3,1)$ or ${\rm SO}(4,1)$.
Denote $\Gamma_\infty={\rm Stab}_\Gamma(p_\infty)$; by assumption $\Gamma_\infty=\la T,U \ra$ with $T,U$ translations in  $\partial_\infty \HR^3 \setminus \{ p_\infty \} \simeq \R^2$. Assume that $T \in \Delta_\infty={\rm Stab}_\Delta(p_\infty)$ and that $U$ is the stable letter in the HNN decomposition of $\Gamma$ given by $\Sigma$.
\begin{enumerate}

\item[(a)]  If $G={\rm SU}(3,1)$ then $\rho_t(U)$ is ellipto-parabolic for all $t \neq 0$ in a neighborhood $I$ of $0$, hence $\rho_t$ is parabolic-preserving for $t \in I$. The restriction of $\rho_t$ to $\Gamma_\infty$ is faithful for all $t \in I$; moreover, if $T,U$ are orthogonal it is also discrete, hence $\rho_t$ is strongly parabolic-preserving.

\item[(b)] If $G={\rm SO}(4,1)$ and $T,U$ are orthogonal, then $\rho_t(U)$ is elliptic for all $t \neq 0$ in a neighborhood of $0$. In particular, $\rho_t$ is not parabolic-preserving for $t \neq 0$ in this neighborhood.

\item[(c)]  If $G={\rm SO}(4,1)$ and $T,U$ are not orthogonal, then $\rho_t(U)$ is ellipto-parabolic for all $t \neq 0$ in a neighborhood of $0$. In this case, $\rho_t$ is parabolic-preserving for $t \neq 0$ in this neighborhood, but not strongly parabolic-preserving.
\end{enumerate}
\end{prop}

Note: The assumption that $U$ is the stable letter in the HNN decomposition of $\Gamma = \pi_1(M)$ given by $\Sigma$ means geometrically that the horospherical--geodesic loop representing $U \in \pi_1(M)$ intersects $\Sigma$ only once (see Lemma~\ref{CuspStableLetter} below). While this assumption together with the assumption that $T,U$ are orthogonal with $T \in \Delta_\infty$ may seem very restrictive, we will see in Section~\ref{BianchiBending1} that they apply to the Bianchi groups, providing infinitely many examples in dimension 3, even up to commensurability.

\medskip

\Pf In the Siegel model for $\HR^3$ (taking $p_\infty = (1,0,0,0)^T$) we may normalize the translations $T,U$ to have the form:

\begin{equation}\label{CuspH3RGenerators}
\begin{array}{lr}
T=\left( \begin{array}{cccc}
1 & -a & 0 & -a^2/2\\
0 & 1 & 0 & a \\
0 & 0 & 1 & 0 \\
0 & 0 & 0 & 1 
\end{array}\right) ,
&
U=\left( \begin{array}{cccc}
1 & -b_1 & -b_2 & -(b_1^2+b_2^2)/2 \\
0 & 1 & 0 & b_1 \\
0 & 0 & 1 & b_2 \\
0 & 0 & 0 & 1 
\end{array}\right) 
\end{array}
\end{equation}

for some $a,b_1,b_2 \in \R$ with $a \neq 0, b_2 \neq 0$. Note that $T,U$ are orthogonal if and only if $b_1=0$. 

Moreover, denoting $(e_1,e_2,e_3,e_4)$ the standard basis of $\R^4$, $T$ preserves $e_3^\perp$, whose image in $\HR^3$ is a totally geodesic copy of $\HR^2$, and we may assume that this is the copy of $\HR^2$ preserved by the surface subgroup $\Delta$. 

Now the centralizer of $\Delta$ in $G={\rm PU(3,1)}$ is the 1-parameter group $Z(\Delta)=\{ Z_3(u) \, | \, u \in {\rm U}(1)\}$, where $Z_3(u)$ denotes the diagonal matrix ${\rm Diag}(1,1,u,1)$. (For simplicity we do not normalize this matrix to lie in SU(3,1)).

Likewise, embedding (the Siegel model of) SO(3,1) into SO(4,1) by embedding $\R^{3,1}$ into $\R^{4,1}$ as ${\rm Span}(e_1,e_2,e_3,e_5)$,  the centralizer of $\Delta$ in $G={\rm SO(4,1)}$ is the 1-parameter group $Z(\Delta)=\{ R_{34}(t) \, | \, t \in \R \}$ where, denoting for convenience $c_t={\rm cos} \, t$ and $s_t={\rm sin} \, t$:
$$
R_{34}(t)=\left( \begin{array}{ccccc}
1 & 0 & 0 & 0 & 0 \\
0 & 1 & 0 & 0 & 0 \\
0 & 0 & c_t & -s_t & 0\\
0 & 0 & s_t & c_t & 0 \\
0 & 0 & 0 & 0 & 1 
\end{array}\right)
$$

\begin{enumerate}

\item[(a-b-c)] {\bf Parabolicity:}  By (\ref{BendingHNN}), the bending deformations $\rho_t$ act on $\Gamma_\infty=\la T,U \ra$ by $\rho_t(T)=T$ and $\rho_t(U)=g_tU$, where $(g_t)_{t \in \R}$ is a smooth parametrization of $Z(\Delta)$ as above. More explicitly, in the case where $G={\rm SU}(3,1)$, denoting $u=e^{it} \in {\rm U(}1)$ and $Z_3(u)={\rm Diag}(1,1,u,1)$ as above:

$$
\rho_t(U)=Z_3(u)U=\left( \begin{array}{cccc}
1 & -b_1 & -b_2 & -(b_1^2+b_2^2)/2 \\
0 & 1 & 0 & b_1 \\
0 & 0 & u & ub_2 \\
0 & 0 & 0 & 1 
\end{array}\right) 
$$
(As above, for simplicity we do not normalize this matrix to lie in ${\rm SU}(3,1)$). 
This is an ellipto-parabolic isometry of ${\rm H}^3_\C$ whenever $u \neq 1$. Indeed, it acts on $\partial_\infty {\rm H}^3_\C \setminus \{ \infty \} \simeq \C^2 \times \R$ in Heisenberg coordinates by:

$$ (z_1,z_2,v) \longmapsto (z_1+b_1,u(z_2+b_2),v+2b_1 {\rm Im} (z_1)+2b_2 {\rm Im}(z_2))
$$
This has no fixed point other than $\infty$ in $\partial_\infty {\rm H}^3_\C$ when $u \neq 1$. Indeed, if $b_1 \neq 0$ the first coordinate is never fixed.
If $b_1=0$, the second coordinate is fixed exactly when $z_2=\frac{ub_2}{1-u}$. But this has ${\rm Im}(z_2)\neq 0$ so the third coordinate is not fixed, as $b_2 \neq 0$.

(We will discuss discreteness in the next item).

Likewise, in the case where $G={\rm SO}(4,1)$, with the above notation:

$$
\rho_t(U)=R_{34}(t)U=\left( \begin{array}{ccccc}
1 & -b_1 & -b_2 & 0 & -(b_1^2+b_2^2)/2 \\
0 & 1 & 0 & 0 & b_1 \\
0 & 0 & c_t & -s_t &c_t b_2 \\
0 & 0 & s_t & c_t & s_t b_2 \\
0 & 0 & 0 & 0 & 1 
\end{array}\right)
$$
In contrast with the previous case, when $u \neq 1$ this is an elliptic isometry of ${\rm H}^4_\R$ if $b_1=0$, otherwise it is ellipto-parabolic. Indeed, note that it acts on $\partial_\infty {\rm H}^4_\R \setminus \{ \infty \} \simeq \R^3$ in Heisenberg coordinates by:

$$ (x,y,z) \longmapsto (x+b_1,c_ty-s_tz+c_tb_2,s_ty+c_tz+s_tb_2)
$$
This is a rotation (through angle $t$) in the $(y,z)$-plane when $b_1=0$, and a screw-motion of $\R^3$ when $b_1 \neq 0$.

\item[(a)] {\bf Discreteness/faithfulness:}  As noted above, $\rho_u(T)=T$ and $\rho_u(U)$ act on $\partial_\infty {\rm H}^3_\C \setminus \{ \infty \} \simeq \C^2 \times \R$ in Heisenberg coordinates by:

$$
\begin{array}{c} T: (z_1,z_2,v) \longmapsto (z_1+a,z_2,v+2a {\rm Im} (z_1)) 
\\
\\
\rho_u(U):  (z_1,z_2,v) \longmapsto (z_1+b_1,u(z_2+b_2),v+2b_1 {\rm Im} (z_1)+2b_2 {\rm Im}(z_2))
\end{array}
$$
Note that the action of $\rho_u(U)$ on the $z_2$ coordinate is (for $u \neq 1$) a rotation $z_2 \mapsto u(z_2+b_2)$ with angle $t$ and center $c_u=\frac{ub_2}{1-u}$. Note that ${\rm Im} (c_u) \neq 0$ for $u \neq -1$. Changing the $z_2$ coordinate to $z'_2=z_2-c_u$, the action of $\rho_u(U)$ is now given by: 

$$\rho_u(U):  (z_1,z'_2,v) \longmapsto (z_1+b_1,uz'_2,v+2b_1 {\rm Im} (z_1)+2b_2 {\rm Im}(z'_2)+2b_2 {\rm Im}(c_u))
$$
 
\begin{lem}\label{SU(3,1)StronglyParabolic} For any $u \neq 1 \in {\rm U}(1)$, the restriction of $\rho_u$ to $\Gamma_\infty = \la T,U \ra$ is faithful. When $T,U$ are orthogonal it is also discrete.
\end{lem} 

\Pf (of Lemma~\ref {SU(3,1)StronglyParabolic}) Using the fact that $\rho_u(T)=T$ and $\rho_u(U)$ commute, the orbit of any given point  $p_0=(z_1,z'_2,v) \in \partial_\infty {\rm H}^3_\C \setminus \{ \infty \}$ is $\{T^m\rho_u(U)^n p_0\ \, | \, m,n \in \Z \}$, where:

\begin{equation}\label{orbit}T^m\rho_u(U)^n p_0 = \big(z_1+ma+nb_1,u^nz'_2,v+2(ma+nb_1) {\rm Im} (z_1)+2nb_2 {\rm Im}(c_u)+2b_2\sum_{j=0}^{n-1}  {\rm Im}(u^jz'_2) \big).
\end{equation}

In particular, taking $p_0=(0,0,0)$ gives: $T^m\rho_u(U)^n p_0 = \big(ma+nb_1,0,2nb_2 {\rm Im}(c_u))$. As noted above, ${\rm Im}(c_u) \neq 0$ (and $a \neq 0, b_2 \neq 0$), hence if $T^m\rho_u(U)^n p_0 =p_0$, then $n=m=0$. This proves faithfulness of $\rho_u$ restricted to $\Gamma_\infty$.

We now show discreteness of $\rho_u(\Gamma_\infty)$ assuming that $T,U$ are orthogonal, that is when $b_1=0$.
More specifically we show that the orbits in $\partial_\infty {\rm H}^3_\C \setminus \{ \infty \}$ are all discrete,

Assume that $T^{m_k}U^{n_k}(p_0)$ converges to some $q=(w_1,w_2,s)$ in the above coordinates, along some sequence $(m_k, n_k)_{k\in\N}$ in $\Z^2$ diverging to $\infty$. In view of (\ref{orbit}), this means that, as $k \longrightarrow \infty$:

$$\left\lbrace\begin{array}{l}
(1) \ z_1+m_ka+n_kb_1 \longrightarrow w_1 \\
(2) \ u^{n_k}z'_2  \longrightarrow w_2 \\
(3) \ v+2(m_ka+n_kb_1) {\rm Im} (z_1)+2n_kb_2 {\rm Im}(c_u)+2b_2\sum_{j=0}^{k-1}  {\rm Im}(u^{n_j}z'_2) \longrightarrow s
\end{array}\right.
$$

Now by assumption $b_1=0$, so condition (1) reduces to: $\ z_1+m_ka \longrightarrow w_1$, which implies that $m_k$ is eventually constant, say $m_k=M$ for all $k$ large enough. But then $T^{m_k}U^{n_k}(p_0)=T^MU^{n_k}(p_0)=U^{n_k}(T^M(p_0))$ is contained in an orbit of the cyclic subgroup $\la U \ra$, which is discrete as $U$ has infinite order.
 
\item[(c)] {\bf Discreteness/faithfulness:} Now assume that $G={\rm SO}(4,1)$ and $T,U$ are not orthogonal, that is $b_1 \neq 0$.

Recall from above that $\rho_t(T)$ acts on $\partial_\infty {\rm H}^4_\R \setminus \{ \infty \} \simeq \R^3$ in Heisenberg coordinates as the translation $ (x,y,z) \longmapsto (x+a, y, z)$, and $\rho_t(U)$ acts as the screw-motion:
$$ (x,y,z) \longmapsto (x+b_1,c_ty-s_tz+c_tb_2,s_ty+c_tz+s_tb_2)
$$
By Lemma~\ref{RtimesS1} $\rho_t$ is a non-discrete or non-faithful representation of $\Gamma_\infty \simeq \Z^2$.
\end{enumerate}
\EPf

The same result holds in any dimension $n \geqslant 3$, with the same proof (note that the relevant centralizer $Z(\Delta)$ is still 1-dimensional, hence the effect of bending on matrices generating the cusp group is still restricted to the corresponding one or two lines as above).  

\begin{prop}\label{CuspBendingNonSeparatingDimn} Let $M=\HR^n/\Gamma$ ($n \geqslant 3$) be a finite-volume hyperbolic n-orbifold with a torus cusp represented by $p_\infty \in \partial_\infty \HR^n$, and $\Sigma=\HR^{n-1}/\Delta$ a properly embedded, non-separating totally geodesic hypersurface in $M$ containing (the image of) $p_\infty$. 
Let $\rho_t: \Gamma \longrightarrow G$ ($t \in \R$) denote the bending deformations of $\Gamma$ along $\Delta$, with $G={\rm SU}(n,1)$ or ${\rm SO}(n+1,1)$.
Denote $\Gamma_\infty={\rm Stab}_\Gamma(p_\infty)$; by assumption $\Gamma_\infty=\la T_1,...,T_{n-2},U \ra$ with $T_i,U$ translations in  $\partial_\infty \HR^n \setminus \{ p_\infty \} \simeq \R^{n-1}$. Assume that $\Delta_\infty={\rm Stab}_\Delta(p_\infty)=\la T_1,...,T_{n-2} \ra$ and that $U$ is the stable letter in the HNN decomposition of $\Gamma$ given by $\Sigma$.
\begin{enumerate}

\item[(a)]  If $G={\rm SU}(n,1)$ then $\rho_t(U)$ is ellipto-parabolic for all $t \neq 0$ in a neighborhood $I$ of $0$, hence $\rho_t$ is parabolic-preserving for $t \in I$. The restriction of $\rho_t$ to $\Gamma_\infty$ is faithful for all $t \in I$; moreover, if $U$ is orthogonal to ${\rm Span} ( T_1,...,T_{n-2})$ it is also discrete, hence $\rho_t$ is strongly parabolic-preserving.

\item[(b)] If $G={\rm SO}(n+1,1)$ and $U$ is orthogonal to ${\rm Span} ( T_1,...,T_{n-2})$, then $\rho_t(U)$ is elliptic for all $t \neq 0$ in a neighborhood of $0$. In particular, $\rho_t$ is not parabolic-preserving for $t \neq 0$ in this neighborhood.

\item[(c)]  If $G={\rm SO}(n+1,1)$ and $U$ is not orthogonal to ${\rm Span} ( T_1,...,T_{n-2})$,  then $\rho_t(U)$ is ellipto-parabolic for all $t \neq 0$ in a neighborhood of $0$. In this case, $\rho_t$ is parabolic-preserving for $t \neq 0$ in this neighborhood, but not strongly parabolic-preserving.
\end{enumerate}
\end{prop}

Note: As above, the assumption that $U$ is the stable letter in the HNN decomposition of $\Gamma = \pi_1(M)$ given by $\Sigma$ may seem very restrictive, but we will see in Section~\ref{BendingAllDims} that it applies to infinitely many pairwise non-commensurable groups in any dimension.

\begin{lem}\label{CuspStableLetter} With the above notation, if $\Sigma$ intersects a cusp cross-section of $M$ in a single component, then ($\Sigma$ is non-separating and) one can choose generators $(T_1,...,T_{n-2},U)$ of $\Gamma_\infty$ such that 
$T_1,...,T_{n-2} \in \Delta_\infty$ and $U$ is the stable letter in the HNN decomposition of $\Gamma$ given by $\Sigma$.
\end{lem}

Note: Since $\Sigma$ is properly embedded in $M$, its intersection with an (embedded) $(n-1)$-torus cross-section of $M$ is a disjoint union of parallel $(n-2)$-tori, see e.g. Section 6 of \cite{BM}. The assumption here is that this disjoint union is a single $(n-2)$-torus. This implies that $\Sigma$ does not separate $M$ - if it did, it would also separate the torus cross-section, which requires at least two parallel $(n-2)$-tori.

\medskip

\Pf It is helpful to represent elements of $\Gamma= \pi_1(M)$ by loops in $M$. 
Recall that, while loxodromic elements of $\Gamma= \pi_1(M)$ have a unique geodesic representative in their (free) homotopy class, parabolic elements do not have any geodesic representative. However, in any horosphere preserved by such a parabolic element $\gamma \in \Gamma$ there is a loop $c_\gamma$ representing $\gamma$ which is a geodesic in the induced (Euclidean) metric (unique after choice of basepoint). We will call $c_\gamma$ the  \emph{horospherical--geodesic loop representing $\gamma=[c_\gamma]$.}  

Let $\tilde{H} \subset \HR^n$ be a horosphere based at $p_\infty$ such that $H=\tilde{H}/\Gamma_\infty$ is embedded in $M=\HR^n/\Gamma$, and choose a basepoint $x_0 \in \Sigma \cap H$. 
Choose generators $T_1,...,T_{n-2}$ for $\Delta_\infty \simeq \Z^{n-2}$ and $U \in \Gamma_\infty$ such that $T_1,...,T_{n-2},U$ generate $\Gamma_\infty \simeq \Z^{n-1}$, and denote by $c_{U}$ the horospherical--geodesic loop representing $U$ based at $x_0$. Since $\Sigma$ intersects the cusp cross-section of $M$ in a single component, $c_U$ intersects $\Sigma$ exactly once (at $x_0$). 

Since $\Sigma$ is non-separating, $M'=M \setminus \Sigma$ is the interior of a connected manifold $\bar{M'}$ whose boundary consists of two isometric copies $\Sigma_1$ and $\Sigma_2$ of $\Sigma$. The image of the loop $c_U$ in $M'$ is a simple arc connecting $\Sigma_1$ to $\Sigma_2$, whose interior lies in $M'$; by (the non-separating version of) the van Kampen theorem, $\pi_1(M)$ decomposes as an HNN extension $\pi_1(M')*_z$ with $z=[c_U]=U$. 
\EPf

\subsection{Examples in all dimensions}\label{BendingAllDims}

In this section we construct additional examples in arbitrary dimensions, proving Theorem~\ref{df_all_dims} stated in the introduction. We start with a cusped hyperbolic manifold  $M=\HR^n/\Gamma$, and $\rho_{hyp}: \Gamma \longrightarrow {\rm SO}(n,1)$ the hyperbolic representation of $\Gamma=\pi_1(M)$, i.e. the holonomy representation of the complete hyperbolic structure on $M$.

\medskip

\Pf (of Theorem~\ref{df_all_dims}) We proceed by constructing infinitely many commensurability classes of cusped hyperbolic manifolds containing totally geodesic hypersurfaces. This is done via a well-known arithmetic construction (see e.g. \cite{BHW}). The rough idea is to look at the group, $\Gamma$, of integer points of the orthogonal groups of various carefully selected quadratic forms of signature $(n,1)$. The quotient $M=\HR^n/\Gamma$ will be a cusped hyperbolic $n$-orbifold containing a totally geodesic hypersurface. After passing to a carefully selected cover we can  produce our parabolic preserving representations via the bending construction. 

We now discuss the details for a specific form and observe that the proof is essentially unchanged if one selects a different form. Let $\hat \Gamma={\rm SL}(n+1,\Z)\cap {\rm SO}^0(n,1)$ and let $\hat \Delta=\hat \Gamma\cap {\rm SO}^0(n-1,1)$. The group $\hat \Gamma$ clearly contains unipotent elements and so we see that $\hat M=\HR^n/\hat \Gamma$ is a cusped hyperbolic $n$-orbifold, which contains an immersed totally geodesic codimension-1 suborbifold isomorphic to $\hat \Sigma=\HR^{n-1}/\hat \Delta$. By work of McReynolds--Reid--Stover (Theorem 1.2, Remark 1 and Proposition 3.1 of \cite{MRS}, building upon work of Bergeron--Haglund--Wise, \cite{BHW}) we can find finite index subgroups $\Gamma\subset \hat \Gamma$ and $\Delta\subset \hat \Delta$ and corresponding manifolds $M=\HR^n/\Gamma$ and $\Sigma=\HR^{n-1}/\Delta$ with the following properties. 

\begin{itemize}
\item $\Gamma$ is torsion-free,
\item $\Sigma$ is embedded in $M$,
\item $M$ has only torus cusps, and 
\item $\Sigma$ intersects each cusp cross-section of $M$ in at most one component.
\end{itemize}

Each $M$ contains the totally geodesic hypersurface $\Sigma$ along which we can bend to produce a family $\rho_\theta$ of representations from $\Gamma$ into ${\rm SU}(n,1)$ or ${\rm SO}(n+1,1)$. 
By Lemma~\ref{CuspStableLetter} and Proposition~\ref{CuspBendingNonSeparatingDimn}, the representations $\rho_\theta$ into ${\rm SU}(n,1)$ are then parabolic-preserving (while those into ${\rm SO}(n+1,1)$ are not strongly parabolic-preserving). . \EPf

\begin{rmk}
\begin{itemize}
\item[(a)] Proposition~\ref{CuspBendingNonSeparatingDimn} would give us the stronger result that the deformations into ${\rm SU}(n,1)$ are strongly parabolic-preserving, if we could certify the orthogonality condition on every cusp. While this can be easily shown for a single cusp (e.g. the standard cusp for the above lattice ${\rm SL}(n+1,\Z)\cap {\rm SO}^0(n,1)$), it seems difficult when $M$ has more than one cusp, which must be the case in higher dimensions by the main result of \cite{St}. 

\item[(b)] It is well known, see \cite{Th} or more generally \cite{HT}, that the complement in $S^3$ of the figure-8 knot does not contain an embedded totally geodesic hypersurface. Therefore, the deformations produced in Theorem \ref{figure8} are distinct from those produced by Theorem \ref{df_all_dims}.
\end{itemize}
\end{rmk}

\subsection{Bending the Bianchi orbifolds into SU(3,1) along the modular surface}\label{BianchiBending1}

We now compute explicit examples in dimension 3 which illustrate the above construction. 
The Bianchi groups are defined, for $d\geqslant 1$ squarefree, by ${\rm Bi}(d)={\rm PSL}(2,\mathcal{O}_d) < {\rm PSL}(2,\C)$, where $\mathcal{O}_d$ denotes as usual the ring of integers of $\Q (i\sqrt{d})$. 
 These are well known to be arithmetic lattices in  ${\rm PSL}(2,\C)$, with number of cusps equal to the class number $h_d$ of  $\Q (i\sqrt{d})$. In particular, by the Baker--Heegner--Stark theorem, ${\rm Bi}(d)$ has a single cusp exactly when $d=1,2,3,7,11,19,43,67,163$. It is also known that the cusp cross-sections of ${\rm Bi}(d)$ are all tori when $d \neq 1,3$. 

The Bianchi groups can be regarded as subgroups of ${\rm SO}(3,1)$, since ${\rm PSL}(2,\C)$ and ${\rm SO}^0(3,1)$ are isomorphic Lie groups (they are both isomorphic to ${\rm Isom}^+({\rm H}^3_\R)$). One can then ask about the existence and potential properties of deformations into ${\rm SL}(4,\R)$ or ${\rm SU}(3,1)$. These were studied by the second author and Thistlethwiate in \cite{PT} for small values of $d$. In particular, it was shown there that the Bianchi groups ${\rm Bi}(1)$ and ${\rm Bi}(3)$ have no parabolic-preserving deformations into ${\rm SU}(3,1)$. More specifically, ${\rm Bi}(3)$ has a 1-dimensional space of deformations into ${\rm SU}(3,1)$, none of which are parabolic-preserving, whereas any deformation of ${\rm Bi}(1)$ into ${\rm SU}(3,1)$ is conjugate to one into ${\rm SO}(3,1)$, and these were known to not be parabolic-preserving.

The next simplest cases correspond to $d=2,7,11$; in \cite{PT} it was shown that the corresponding deformation spaces of ${\rm Bi}(d)$ into ${\rm SU}(3,1)$ have dimension at most 3, 3, 4 respectively. The bending construction described above and made explicit below shows that these spaces contain a smooth 1-parameter family of deformations transverse to the 2-dimensional family of Dehn surgery deformations into ${\rm SO}(3,1)$. We will give more details for these three cases below, as the generators and presentations are especially simple.

  Swan found explicit presentations for most Bianchi groups ${\rm Bi}(d)$ with $d \leqslant 19$ in \cite{Sw}, from which we review a few basic facts. Recall that $\mathcal{O}_d=\Z[\tau]$, where $\tau=i\sqrt{d}$ if $d=1,2$ mod 4 and $\tau=\frac{1+i\sqrt{d}}{2}$ if $d=3$ mod 4. 
     
   The generators of  ${\rm Bi}(d)$ include for all $d$:
$$\begin{array}{ccc}
T=\left(\begin{array}{cc}
1 & 1 \\
0 & 1 \end{array}\right)
&
U=\left(\begin{array}{cc}
1 & \tau \\
0 & 1 \end{array}\right)
&
A=\left(\begin{array}{cc}
0 & -1 \\
1 & 0 \end{array}\right)
\end{array}
$$

These particular generators are especially relevant for our purposes, since $T,U$ generate the stabilizer of the standard cusp, and $T,A$ generate the modular group ${\rm PSL}(2,\Z)$.
It turns out that these suffice to generate the whole group ${\rm Bi}(d)$ when $d=2,7,11$. The corresponding presentations given by Swan are the following:

\begin{equation}\label{Bianchi2711Pres}
\begin{array}{l}
{\rm Bi}(2) = \langle a, t, u \, | \, [t,u]=a^2=(at)^3=(au^{-1}au)^2=1 \rangle \\
\\
{\rm Bi}(7)= \langle a, t, u \, | \, [t,u]=a^2=(at)^3=(atu^{-1}au)^2=1 \rangle \\
\\
{\rm Bi}(11)= \langle a, t, u \, | \, [t,u]=a^2=(at)^3=(atu^{-1}au)^3=1 \rangle 
\end{array}
\end{equation}

In order to perform explicit computations, we first realize the Bianchi groups (and in particular the above generators) as subgroups/elements of SO(3,1). For this we use an explicit isomorphism ${\rm PSL}(2,\C) \longrightarrow {\rm SO}^0(3,1)$, for example the Spin homomorphism given in \cite{FS}. This lands us in the ball model of ${\rm H}^3_\R$, that is in the group ${\rm SO}(Q_1)$ where $Q_1$ is the quadratic form $Q_1(x,y,z,t)=x^2+y^2+z^2-t^2$. For our purposes it will be more convenient to work in the \emph{Siegel model} corresponding to the quadratic form $Q_2(x,y,z,t)=2xt+y^2+z^2$, where the standard parabolic subgroup is upper-triangular.

In these coordinates, the lattice embedding $\rho_1$ of ${\rm Bi}(d)$ into SO(3,1) is given by $\rho_1(a)=A_1$, $\rho_1(t)=T_1$, and $\rho_1(u)=U_1$ , where:

\begin{equation}\label{BianchiGens1}
\begin{array}{lr}
A_1=\left( \begin{array}{cccc}
0 & 0 & 0 & -1 \\
0 & -1 & 0 & 0 \\
0 & 0 & 1 & 0 \\
-1 & 0 & 0 & 0 
\end{array}\right) 
&
T_1=\left( \begin{array}{cccc}
1 &- \sqrt{2} & 0 & -1 \\
0 & 1 & 0 & \sqrt{2} \\
0 & 0 & 1 & 0 \\
0 & 0 & 0 & 1 
\end{array}\right) 
\end{array}
\end{equation}

\begin{equation}\label{BianchiGens2}
U_1=\left( \begin{array}{cccc}
1 & 0 & \sqrt{2d} & -d \\
0 & 1 & 0 & 0 \\
0 & 0 & 1 &  -\sqrt{2d} \\
0 & 0 & 0 & 1 
\end{array}\right) \ {\rm if} \ d =1,2 \, mod. 4
\end{equation}

\begin{equation}\label{BianchiGens3}
U_1=\left( \begin{array}{cccc}
1 & -\sqrt{2}/2 &  \sqrt{2d}/2 & -(d+1)/2 \\
0 & 1 & 0 & \sqrt{2}/2 \\
0 & 0 & 1 &  -\sqrt{2d}/2 \\
0 & 0 & 0 & 1 
\end{array}\right) \ {\rm if} \ d =3 \, mod. 4
\end{equation}

It is well known that the modular surface ${\rm H}^2_\R/{\rm PSL}(2,\Z)$ is embedded in the Bianchi orbifold ${\rm H}^3_\R/{\rm Bi}(d)$ for $d\neq 1,3$ and is non-separating, see for example Section 6.3.2 of \cite{F}. For $d=2, 7, 11$ Fine also gives in Section 4.4 of \cite{F} the details of the corresponding HNN extension structure. In the notation of the presentation \ref{Bianchi2711Pres} above, with $v=u^{-1}au$ this gives:

$$ {\rm Bi}(d)=\Gamma_1*_u
$$ 
where $\Gamma_1=\langle a,t,v \rangle$ has isomorphic subgroups $H=\langle a,t \rangle={\rm PSL}(2,\Z)$ and $H'=\langle v,t \rangle$, with isomorphism realized in ${\rm Bi}(d)$ by conjugation by $u$. (Note that $v=u^{-1}au$ and $t=u^{-1}tu$). 

Now the centralizer of $H=\langle A_1, T_1 \rangle $ in $G={\rm PU(3,1)}$ is the 1-parameter group $Z(H)=\{ Z_3(u) \, | \, u \in {\rm U}(1)\}$, where $Z_3(u)$ denotes the diagonal matrix ${\rm Diag}(1,1,u,1)$. (For simplicity we do not normalize this matrix to lie in SU(3,1)). We therefore obtain the following 1-parameter family $\rho_u$ ($u \in {\rm U}(1)$) of deformations of the $\R$-Kleinian embedding. 
Define, for $u \in {\rm U}(1)$, $\rho_u$ by $\rho_u(a)=A_1$, $\rho_u(t)=T_1$ and $\rho_u(u)=U_u=Z_3(u)U_1$ where $A_1, T_1$ are as above and:

\begin{equation}\label{BianchiBend1}
U_u=\left( \begin{array}{cccc}
1 & 0 & \sqrt{2d} & -d \\
0 & 1 & 0 & 0 \\
0 & 0 & u &  -u\sqrt{2d} \\
0 & 0 & 0 & 1 
\end{array}\right) \ {\rm if} \ d =1,2 \, mod. 4
\end{equation}

\begin{equation}\label{BianchiBend2}
U_u=\left( \begin{array}{cccc}
1 &- \sqrt{2}/2 &  \sqrt{2d}/2 & -(d+1)/4 \\
0 & 1 & 0 & \sqrt{2}/2 \\
0 & 0 & u &  -u\sqrt{2d}/2 \\
0 & 0 & 0 & 1 
\end{array}\right) \ {\rm if} \ d =3 \, mod. 4
\end{equation}

Note that, by construction, the modular subgroup $\rho_u(\langle a, t \rangle)=\langle A_1,T_1\rangle$ remains constant for all $u \in {\rm U}(1)$. Also, as previously, for simplicity we do not normalize these matrices to lie in ${\rm SU}(3,1)$. 

We summarize the relevant properties of these bending representations below. Parts 1 and 2 follow from the results of \cite{JM}, with the caveat that the Bianchi orbifolds are orbifolds and not manifolds -- hence we give an independent proof.

\begin{prop}\label{BianchiDefsSU31}
\begin{enumerate}

\item For any $u \in {\rm U}(1)$ and $d=2,7,11$, the matrices $A_1$, $T_1$ and $U_u$ defined above define a representation $\rho_u: {\rm Bi}(d) \longrightarrow {\rm U}(3,1)$. For any $d \neq 1,3$ they extend to a representation $\rho_u: {\rm Bi}(d) \longrightarrow {\rm U}(3,1)$.

\item For any $d \neq 1,3$ and any $u \neq u' \in {\rm U}(1)$, the representations $\rho_u$ and $\rho_{u'}$ are not conjugate.

\item For any $d \neq 1,3$ and any $u \in {\rm U}(1)$, $\rho_u$ is parabolic-preserving. 

\item For any $d \neq 1$ congruent to 1 or 2 mod. 4, and any $u \in {\rm U}(1)$, $\rho_u$ is strongly parabolic-preserving. 

\item For any $d \neq 1,3$ and any $u \neq \pm 1 \in {\rm U}(1)$, $\rho_u({\rm Bi}(d))$ is not contained in any proper Lie subgroup of ${\rm U}(3,1)$, and in particular is Zariski-dense in ${\rm U}(3,1)$.

\end{enumerate}
\end{prop}

\Pf \begin{enumerate}
\item This follows from the HNN decomposition described above, but a straightforward matrix computation, for example using formal algebra software, also shows that these matrices satisfy the relations in the above presentations of ${\rm Bi}(d)$ when $d=2,7,11$.

\item This follows from the observation that ${\rm Tr} \, U_u=3+u$ for all $u \in {\rm U}(1)$.

\item The representation $\rho_u$ is parabolic-preserving because $U_u$ is not diagonalizable, hence parabolic, as is easily seen in this upper-triangular form.

\item When $d \neq 1$ is congruent to 1 or 2 mod. 4, the above generators $T_1,U_1$ of the cusp group are orthogonal, hence $\rho_u$ is strongly parabolic-preserving by Proposition~\ref{CuspBendingNonSeparatingDim3} (a).

\item Let $G$ be the smallest connected Lie subgroup of ${\rm U}(3,1)$ containing $\rho_u({\rm Bi}(d))$. By Theorem 4.4.2 of \cite{CG}, either $G={\rm U}(3,1)$ or $G$ has a global fixed point in $\overline{{\rm H}_\C^3}$ or preserves a totally geodesic submanifold of ${\rm H}_\C^3$. By the classification of totally geodesic subspaces (Proposition 2.5.1 of \cite{CG}), in all these cases $G$ would be reducible, or contained in a conjugate of ${\rm O}(3,1)$. It is easy to see that $\rho_u$ is irreducible for all $u \in {\rm U}(1)$. Moreover, since ${\rm Tr} \, U_u=3+u$, $\rho_u({\rm Bi}(2))$ is not contained in a conjugate of ${\rm O}(3,1)$ unless $u=\pm 1$. 
\end{enumerate}
\EPf

 \subsection{Bending the Bianchi orbifolds into SO(4,1) along the modular surface}\label{BianchiBending2}
 
 We embed (the Siegel model of) SO(3,1) into SO(4,1) by embedding $\R^{3,1}$ into $\R^{4,1}$ as ${\rm Span}(e_1,e_2,e_3,e_5)$.
The elements $A_1,T_1,U_1$ of ${\rm Bi}(d)$ in SO(3,1) given in (\ref{BianchiGens1}),(\ref{BianchiGens2}) and (\ref{BianchiGens3}) now become (changing the subscript 1 to 0 as the deformations will now be parametrized by an angle rather than a unit complex number):   
   
   \begin{equation}
\begin{array}{lr}
A'_0=\left( \begin{array}{ccccc}
0 & 0 & 0 & 0 & -1 \\
0 & -1 & 0 & 0 & 0\\
0 & 0 & 1 & 0 & 0\\
0 & 0 & 0 & 1 & 0 \\
-1 & 0 & 0 & 0 & 0
\end{array}\right) 
&
T'_0=\left( \begin{array}{ccccc}
1 & -\sqrt{2} & 0 & 0 & -1 \\
0 & 1 & 0 & 0 & \sqrt{2} \\
0 & 0 & 1 & 0 & 0\\
0 & 0 & 0 & 1 & 0 \\
0 & 0 & 0 & 0 & 1 
\end{array}\right) 
\end{array}
\end{equation}

\begin{equation}
U'_0=\left( \begin{array}{ccccc}
1 & 0 & \sqrt{2d} & 0 & -d \\
0 & 1 & 0 & 0& 0 \\
0 & 0 & 1 &  0 & -\sqrt{2d} \\
0 & 0 & 0 & 1 & 0 \\
0 & 0 & 0 & 0 & 1 
\end{array}\right) \ {\rm if} \ d =1,2 \, mod. 4
\end{equation}

\begin{equation}
U'_0=\left( \begin{array}{ccccc}
1 & -\sqrt{2}/2 &  \sqrt{2d}/2 &0 &  -(d+1)/2 \\
0 & 1 & 0 & 0 & \sqrt{2}/2 \\
0 & 0 & 1 &  0 &-\sqrt{2d}/2 \\
0 & 0 & 0 & 1 & 0 \\
0 & 0 & 0 & 0 & 1 
\end{array}\right) \ {\rm if} \ d =3 \, mod. 4
\end{equation}

Now the centralizer of the modular group $H=\langle A'_0, T'_0 \rangle $ in $G={\rm SO(4,1)}$ is the 1-parameter group $Z(H)=\{ R_{34}(\theta) \, | \, \theta \in \R/2\pi\Z \}$ where, denoting for convenience $c_\theta={\rm cos} \, \theta$ and $s_\theta={\rm sin} \, \theta$:
$$
R_{34}(\theta)=\left( \begin{array}{ccccc}
1 & 0 & 0 & 0 & 0 \\
0 & 1 & 0 & 0 & 0 \\
0 & 0 & c_\theta & -s_\theta & 0\\
0 & 0 & s_\theta & c_\theta & 0 \\
0 & 0 & 0 & 0 & 1 
\end{array}\right)
$$

We therefore obtain the following 1-parameter family $\rho_\theta$ ($\theta \in \R/2\pi\Z$) of deformations into SO(4,1).
Define $\rho_\theta$ by $\rho_\theta(a)=A'_0$, $\rho_\theta(t)=T'_0$ and $\rho_\theta(u)=U'_\theta=R_{34}(\theta)U'_0$ where $A'_0, T'_0$ are as above and:

\begin{equation}
U'_\theta=\left( \begin{array}{ccccc}
1 & 0 & \sqrt{2d} & 0 & -d \\
0 & 1 & 0 & 0& 0 \\
0 & 0 &  c_\theta &   -s_\theta & - c_\theta \sqrt{2d} \\
0 & 0 &  s_\theta &  c_\theta & - s_\theta \sqrt{2d} \\
0 & 0 & 0 & 0 & 1 
\end{array}\right) \ {\rm if} \ d =1,2 \, mod. 4
\end{equation}

\begin{equation}
U'_\theta=\left( \begin{array}{ccccc}
1 & -\sqrt{2}/2 &  \sqrt{2d}/2 &0 &  -(d+1)/2 \\
0 & 1 & 0 & 0 & \sqrt{2}/2 \\
0 & 0 &  c_\theta  &  -s_\theta  &- c_\theta \sqrt{2d}/2 \\
0 & 0 &  s_\theta  &  c_\theta  & - s_\theta \sqrt{2d}/2 \\
0 & 0 & 0 & 0 & 1 
\end{array}\right) \ {\rm if} \ d =3 \, mod. 4
\end{equation}

As in the previous section this defines/extends to a 1-parameter family of representations $\rho_\theta: {\rm Bi}(d) \longrightarrow {\rm SO}(4,1)$. Now however, by Proposition~\ref{CuspBendingNonSeparatingDim3} (b) and (c) we have:

\begin{prop}\label{BianchiDefsSO41} The representations $\rho_\theta: {\rm Bi}(d) \longrightarrow {\rm SO}(4,1)$ are not strongly parabolic-preserving unless $\theta = 0$ or $\pi$.
\end{prop}

\raggedright

\begin{flushleft}
  \textsc{Samuel Ballas\\
   Department of Mathematics, Florida State University}\\
   \verb|sballas@fsu.edu|
\end{flushleft}

\begin{flushleft}
  \textsc{Julien Paupert\\
   School of Mathematical and Statistical Sciences, Arizona State University}\\
       \verb|paupert@asu.edu|
\end{flushleft}

\begin{flushleft}
  \textsc{Pierre Will\\
   Institut Fourier, Universit\'e Grenoble-Alpes}\\
   \verb|pierre.will@univ-grenoble-alpes.fr|
\end{flushleft}

\end{document}